\documentclass[12pt]{article}
\usepackage[T2A]{fontenc}
\usepackage[cp1251]{inputenc}
\usepackage[english,russian]{babel}
\usepackage{amssymb}
\usepackage{amsmath}
\usepackage{amsthm}
\usepackage{amsfonts}
\usepackage{amssymb}

\begin{document}
\language=0
\begin{center}
{\bf TWO PROBLEMS FOR ONE HYPERBOLIC EQUATION\\ OF THE THIRD ORDER IN THREE-DIMENSIONAL SPACE}\\[4mm]        %Заголовок статьи

{\bf Dolgopolov Vyacheslav Mikhailovitch},\\ Department of Mathematics \& Business Informatics,
Samara
State University,

Russia, 443011, Samara,
%%Acad.
Academician Pavlov str.,~1.
%%{V.M.\,Dolgopolov}

{\bf Dolgopolov Mikhail
Vyacheslavovitch} ~
{mvdolg@ssu.samara.ru},\\ Department of General and
Theoretical Physics, Samara State University.
%%{M.V.\,Dolgopolov}

{\bf Irina Nikolayevna
Rodionova},\\ Department of Mathematics \& Business Informatics.
%%{I.N.\,Rodionova}{}
\end{center}

%%\begin{abstract}

\begin{center}
{\sf Abstract}
\end{center}

In the present article, a modified Cauchy problem (problem~$C$) for
the hyperbolic equation of the third order with the data on the
equation's coefficients singularity plane is solved by Riemann method.
The special class in which the solution of the problem $C$ has more
simple appearance is introduced and the area of values of the parameter
$p$ entering into the equation is considerably expanded. In the
special class the mixed problem, which decision was been reduced to the
two-dimensional Volterra's integral equations of the first order with
uncurtailed operators, is considered. Authors found the unique
solution of these equations at various values of the parameter $p$.
\smallskip

This work has been performed within the framework of the project ADTP №\,3341 and 10854

%%\end{abstract}

\smallskip

\noindent {Keywords:} {\sl Volterra integral equations, boundary value problems, hyperbolic type equation.}

\bigskip

\section{Introductory notes}

In the present paper, the following equation is considered:
\begin{equation}
L(U) = U_{xyz} -
\displaystyle\frac{p}{x-y-z}U_{xz}+\displaystyle\frac{p}{x-y-z}U_{yz}-\lambda
U_z = 0
\label{dolgopoloveq:mz1}
\end{equation}
($p$, $\lambda$ "--- const) in the region ${\cal H}=\{(x,y,z):0<z<x-y,
0<y<x<+\infty\}$ of 3-dimensional Euclid space. Using Riemann method,
for the equation (\ref{dolgopoloveq:mz1}), the solution of modified Cauchy problem has been obtained
for the case $0<p<\displaystyle\frac{1}{2}$ (Problem~$C$).
Special representation has been introduced -- special class $W_p$ -- of the solutions of the problem $C$ for equation (\ref{dolgopoloveq:mz1}) with the purpose:\\
1) to simplify the solution form of the problem $C$, making it more convenient to solve new boundary value problems;\\
2) expand the solution of the problem $C$ to the case of negative values of the parameter for which Riemann method is not efficient.

The authors have
obtained such solution of the mixed boundary
value problem for various values of the paameter $p$.

Note that the plane analogue of the equation (\ref{dolgopoloveq:mz1}) with corresponding problems formulation has been considered in the paper~\cite{dolgopolovizv1}.

\section{Problem $\boldsymbol C$}

 In the region ${\cal H}$, find solution of the equation
(\ref{dolgopoloveq:mz1}), continuous in $\overline {\cal H}$, satisfying the conditions:
\begin{equation}
U(x,y,x-y) = \tau(x,y), \qquad 0\leq y\leq x<+\infty;
\label{dolgopoloveq:2}
\end{equation}
\begin{equation}
\lim_{z\to x-y-0} \displaystyle\frac{\partial U}{\partial z} =
\nu(x,y), \qquad 0< y< x<+\infty;
\label{dolgopoloveq:3}
\end{equation}
\begin{equation} \lim_{z\to
x-y-0}(x-y-z)^{2p}(U_{xz}-U_{yz})=\mu(x,y),\quad0<y<x<+\infty\quad(0<p<\displaystyle\frac{1}{2}).
\label{dolgopoloveq:4}
\end{equation}

To solve the problem $C$, we apply Riemann method. Riemann function for the equation~(\ref{dolgopoloveq:mz1})
has been built in the paper \cite{dolgopolovvolko95}. Represent it in the form:
\begin{equation}
\begin{array}{c}
\!\!\!V(x,y,z;x_0,y_0,z_0) =
\displaystyle\frac{(x-y-z)^{2p}}{(x-y_0-z)^{p}(x_0-y-z)^{p}}\times\qquad\qquad\qquad\qquad\qquad\qquad\qquad\\
\quad\times
\displaystyle\sum_{m=0}^\infty\displaystyle\frac{[\lambda(x_0-x)(y_0-y)]^m}{(1)_m
m!}
F\left(p,p+m,1+m;\displaystyle\frac{(x-x_0)(y_0-y)}{(x-y_0-z)(x_0-y-z)}\right),
\end{array}
\label{dolgopoloveq:fRimana}
\end{equation}
where
$$
F(a,b,c,z) =
\displaystyle\sum_{n=0}^\infty\displaystyle\frac{(a)_n(b)_n}{(c)_n
n!} z^n.
$$
In the region ${\cal H}$, take arbitrary point $M_0(x_0,y_0,z_0)$
and, because the coefficients of the equation (\ref{dolgopoloveq:mz1})
\begin{equation} a=
\displaystyle\frac{-p}{x-y-z}, \qquad b =
\displaystyle\frac{p}{x-y-z}
\label{dolgopoloveq:coeffab}\end{equation} become infinity on the
plane $z=x-y$, let us consider the region $H_\varepsilon$ limited by the
planes $x=x_0$, $y=y_0$, $z=z_0$, $z=x-y-\varepsilon$
$(\varepsilon>0)$.

Supposing that the solution of the problem $C$ exists, integrate
Green identity derived in paper \cite{dolgopolovvolko95},
\begin{equation}
V L(U) - U L^*(v) =
\displaystyle\frac{1}{3}\left[\displaystyle\frac{\partial
P}{\partial x} + \displaystyle\frac{\partial Q}{\partial x} +
\displaystyle\frac{\partial H}{\partial x}\right]
\end{equation}
over the region $H_\varepsilon$. For the equation (\ref{dolgopoloveq:mz1}), we have
\begin{equation}
P = V U_{yz} + U V_{yz} - V_yU_z + 3U_z a V,
\label{dolgopoloveq:P8}
\end{equation}
\begin{equation} Q = V(U_{xz} - 3Ub_z) + U V_{xz} - V_zU_x - 3bUV_z,
\label{dolgopoloveq:Q9}
\end{equation}
\begin{equation}
H = V(U_{xy}-3a_xU+3bU_y) - V_xU_y -3aUV_x + V_{xy}U - 3\lambda
U,
\label{dolgopoloveq:H10}
\end{equation}
where $V$ is Riemann function (\ref{dolgopoloveq:fRimana}); $U(x,y,z)$
\cdash--- solution of the equation (\ref{dolgopoloveq:mz1}); $L^*(v)$~--
conjugated operator; coefficients $a$ and $b$ are defined by the formula
(\ref{dolgopoloveq:coeffab}).

Applying Gauss--Ostrogradski
formula to the obtained integral identity, we obtain
\begin{equation}
\displaystyle\sum_{i=0}^3
\displaystyle\int\int\limits_{D_i}
(P\cos\alpha+Q\cos\beta+H\cos\gamma) dS = \displaystyle\sum_{k=1}^4
J_k = 0.
\label{dolgopoloveq:mz11}
\end{equation}
Here, $D_i$ ($i=\overline{0, 3}$)~-- the side of the pyramid
$H_\varepsilon$ located respectively in the planes $x=x_0$,
$y=y_0$, $z=z_0$, $z=x-y-\varepsilon$. Consider each term of the
formula (\ref{dolgopoloveq:mz11})
\begin{equation}
\begin{array}{c}
J_1 = \displaystyle\int\int\limits_{D_0} P|_{x=x_0}
dydz =\displaystyle\int\limits_{z_0}^{x_0-y_0-\varepsilon} dz
\displaystyle\int\limits_{y_0}^{x_0-z-\varepsilon} (V
U_{yz}+U V_{yz}-\qquad\quad\\
\qquad\qquad\qquad\qquad\qquad\qquad\qquad- V_yU_z + 3U_z
aV)|_{x=x_0} dydz = \displaystyle\sum_{x=1}^4 i_k.
\end{array}
\label{dolgopoloveq:mz12}
\end{equation}
In the first two terms, $i_1$ and $i_2$ of eq.
(\ref{dolgopoloveq:mz12}) we integrate by parts taking into account
the property of Riemann function
\begin{equation}
(V_y = aV)_{x=x_0} = 0,
\end{equation}
hence obtain
\begin{equation}
\begin{array}{c}
\!\!\!J_1 = U(x_0,y_0,z_0) - U(x_0,y_0,x_0-y_0-\varepsilon)
+\displaystyle\int\limits_{z_0}^{x_0-z_0-\varepsilon}U_z(x_0,x_0-z-\varepsilon,z)Vdz+\quad\\
\null\hfill+
\displaystyle\int\limits_{y_0}^{x_0-z_0-\varepsilon}U_z(x_0,y,x_0-y-\varepsilon)V_ydy-
\displaystyle\int\limits_{y_0}^{x_0-z_0-\varepsilon}U(x_0,y,z_0)V_ydy.
\end{array}
\label{dolgopoloveq:14}
\end{equation}
Integrating by parts, and
also recalling the relations obtained by direct computation
\begin{equation}
(V_{xz} -
b_zV - bV_z)_{y=y_0} = 0,
\end{equation}
\begin{equation}
aV_x+a_xV+(bV-V_x)_y+\lambda V = 0,
\label{dolgopoloveq:16}
\end{equation} we obtain the following results:
\begin{equation}
\begin{array}{c}
\!\!\!J_2 = - \displaystyle\int\limits_{z_0}^{x_0-y_0-\varepsilon} dz
\displaystyle\int\limits_{z+y_0+\varepsilon}^{x_0} Q(x,y_0,z) dx =
-\displaystyle\int\limits_{z_0+y+\varepsilon}^{x_0}
U_x(x,y_0,x-y_0-\varepsilon)Vdx+\quad\\
\null\hfill+\displaystyle\int\limits_{z_0+y_0+\varepsilon}^{x_0}
U_x(x,y_0,z_0)Vdx -
2\displaystyle\int\limits_{z_0}^{x_0-y_0-\varepsilon}
U(z+y_0+\varepsilon,y_0,z)V_zdz,
\end{array}
\label{dolgopoloveq:17}
\end{equation}

\begin{equation}
\begin{array}{l}
\!\!\!J_3 = -\displaystyle\int\limits_{y_0}^{x_0-z_0-\varepsilon} dy
\displaystyle\int\limits_{z_0+y+\varepsilon}^{x_0} H(x,y,z_0) dx = -
\displaystyle\int\limits_{y_0}^{x_0-z_0-\varepsilon}
U_y(x_0,y,z_0)Vdy+\qquad\qquad\qquad\\
\quad+\displaystyle\int\limits_{y_0}^{x_0-z_0-\varepsilon}
U_y(z_0+y+\varepsilon,y,z_0)Vdy -
\displaystyle\int\limits_{y_0+z_0+\varepsilon}^{x_0}
U(x,x-z_0-\varepsilon,z_0)V_xdx +\\
\null\hfill+\displaystyle\int\limits_{y_0+z_0+\varepsilon}^{x_0}
U(x,y_0,z_0)V_xdx -
3\displaystyle\int\limits_{y_0+z_0+\varepsilon}^{x_0}
U(x,x-z_0-\varepsilon,z_0)(bV-V_x)dx.
\end{array}
\label{dolgopoloveq:18}
\end{equation}

Into the integral
$J_4=\displaystyle\int\limits_{y_0}^{x_0-z_0-\varepsilon}dy\displaystyle\int\limits_{z_0+y+\varepsilon}^{x_0}(-P+Q+H)\vert_{z=x-y-\varepsilon}dx$
, we replace $P$, $Q$, $H$ with their values respectively from equations
(\ref{dolgopoloveq:P8}), (\ref{dolgopoloveq:Q9}),(\ref{dolgopoloveq:H10}) and perform a series of transformations aimed on removing the terms from the double integral containing $U(x,y,x-y-\varepsilon)$. This is
achieved by integrating by parts recalling eq.
(\ref{dolgopoloveq:16}). As a result, we obtain
\begin{equation}\begin{array}{c}J_4=
\displaystyle\int\limits_{y_0}^{x_0-y_0-\varepsilon}dy\displaystyle\int\limits_{z_0+y+\varepsilon}^{x_0}\left[\displaystyle\frac{3}{2}V(U_{xz}-U_{yz})+3(b-a)VU_z+
\displaystyle\frac{3}{2}(V_y-V_x)U_z\right]_{z=x-y-\varepsilon}dx+\\
+\displaystyle\frac{1}{2}\displaystyle\int\limits_{y_0}^{x_0-z_0-\varepsilon}dy\displaystyle\int\limits_{z_0+y+\varepsilon}^{x_0}Vd_x(U_z)+
\displaystyle\frac{1}{2}\displaystyle\int\limits_{y_0}^{x_0-z_0-\varepsilon}\displaystyle\int\limits_{z_0+y+\varepsilon}^{x_0}Vd_y(U_z)dx+\\
+\displaystyle\frac{1}{2}\displaystyle\int\limits_{y_0}^{z_0+y+\varepsilon}dy\displaystyle\int\limits_{z_0+y+\varepsilon}^{x_0}(U_zV_y+V_xU_z)dx+
\displaystyle\int\limits_{y_0+z_0+\varepsilon}^{x_0}U(x,y_0;x-y_0-\varepsilon)V_xdx\\
-\displaystyle\int\limits_{y_0+z_0+\varepsilon}^{x_0}U(x,x-z_0-\varepsilon;z_0)V_xdx
-\displaystyle\int\limits_{y_0}^{x_0-z_0-\varepsilon}U(x_0,y;x_0-y-\varepsilon)V_ydy\\
+\displaystyle\int\limits_{y_0}^{x_0-z_0-\varepsilon}U(z_0+y+\varepsilon,y;z_0)V_ydy
+3\displaystyle\int\limits_{y_0+z_0+\varepsilon}^{x_0}U(x,x-z_0-\varepsilon,z_0)bVdx-\\%%\end{array}$$
-\, 3\displaystyle\int\limits_{y_0+z_0+\varepsilon}^{x_0}U(x,y_0,x-y_0-\varepsilon)bVdx+
\displaystyle\int\limits_{y_0+z_0+\varepsilon}^{x_0}U_x(x,x-z_0-\varepsilon,z_0)Vdx-\\
\qquad\qquad\qquad\qquad\qquad\qquad\qquad\qquad-\displaystyle\int\limits_{y_0+z_0+\varepsilon}^{x_0}U_x(x,y_0,x-y_0-\varepsilon)Vdx.\end{array}\label{dolgopoloveq:19}\end{equation}

The first term of eq. (\ref{dolgopoloveq:19}) will be left intact, as for the 2nd and 3rd, they will be integrated by parts. As a result, all double
integrals in eq. (\ref{dolgopoloveq:19}), except the first will mutually cancel. Substitution of the obtained result and the data of the equations (\ref{dolgopoloveq:14}),
(\ref{dolgopoloveq:17}), (\ref{dolgopoloveq:18}) into the identity (\ref{dolgopoloveq:mz11}),
cancels the single integrals by integrating by parts. The~identity
(\ref{dolgopoloveq:mz11}) takes the form of
\begin{equation}\begin{array}{l}
\!\!\!\!\!\!3U(x_0,y_0,z_0)-3U(x_0,y_0,x_0-y_0-\varepsilon)+
\displaystyle\frac{3}{2}\displaystyle\int\limits_{z_0}^{x_0-y_0-\varepsilon}
U_z(x_0,x_0-z-\varepsilon,z)
\displaystyle\frac{\varepsilon^p}{(x_0-y_0-z)^p} dz+\\
+\displaystyle\frac{3}{2}\displaystyle\int\limits_{z_0}^{x_0-y_0-\varepsilon}U_z(z_0+y_0+\varepsilon,y_0,z)
\displaystyle\frac{\varepsilon^p}{(x_0-y_0-z)^p}dz+\\
\!\!\!\!\!\!+
\displaystyle\frac{3}{2}\displaystyle\int\limits_{y_0}^{x_0-z_0-\varepsilon}dy\displaystyle\int\limits_{z_0+y+\varepsilon}^{x_0}V(U_{xz}-
U_{yz})\vert_{z=x-y-\varepsilon}dx+\qquad\qquad\qquad\qquad\qquad\\
+3\displaystyle\int\limits_{y_0}^{x_0-z_0-\varepsilon}dy\displaystyle\int\limits_{z_0+y+\varepsilon}^{x_0}
[(b-a)V+\displaystyle\frac{1}{2}(V_{y}-V_{x})]U_z\vert_{z=x-y-\varepsilon}dx
=0.\end{array}\label{dolgopoloveq:20}\end{equation}

Performing passage to the limit in the eq. (\ref{dolgopoloveq:20}) at
$\varepsilon\to0$ and recalling the conditions (\ref{dolgopoloveq:2})\,--\,(\ref{dolgopoloveq:4})
we obtain
\begin{equation}
\begin{array}{c}
\!\!\!\!\!\!\!\!\!\!U(x_0,y_0,z_0) = \tau(x_0,y_0) -
\displaystyle\frac{\Gamma(1-2p)}{\Gamma^2(1-p)}
\displaystyle\int\limits_{y_0}^{x_0-z_0}ds\displaystyle\int\limits_{z_0+s}^{x_0}\mu(t,s)\times\qquad\qquad\quad\qquad\\
\qquad\qquad\qquad\times (x_0-t)^{-p}(s-y_0)^{-p} \,{
}_0F_1(1-p;\lambda(x_0-t)(y_0-s))dt -\\
-\displaystyle\frac{\Gamma(2p)}{\Gamma^2(p)}
\displaystyle\int\limits_{y_0}^{x_0-z_0}ds\displaystyle\int\limits_{z_0+s}^{x_0}\nu(t,s)
(x_0-t)^{p-1}(s-y_0)^{p-1} \times\qquad\qquad\qquad\qquad\\
\qquad\qquad\qquad\qquad\times (x_0-y_0-t+s)^{1-2p}\,{
}_0F_1(p;\lambda(x_0-t)(y_0-s))dt,
\end{array}
\label{dolgopoloveq:21}
\end{equation}
$$
{ }_0F_1(\alpha,z) = \displaystyle\sum_{n=0}^\infty
\displaystyle\frac{z^n}{(\alpha)_nn!}.
$$

Designate $D=\{(x,y)/0<y<x<+\infty\}$. Direct checkingmakes us certain that if
\begin{equation}
\tau_{xy} \in C(\overline D); \quad \nu \in C^{(2)}(\overline D); \quad \mu \in C^{(2)}(\overline D),
\label{dolgopoloveq:22}
\end{equation}
then the function (\ref{dolgopoloveq:21}) satisfies the equation
(\ref{dolgopoloveq:mz1}) and conditions (\ref{dolgopoloveq:2}) -- (\ref{dolgopoloveq:4}).

{\bf Theorem 1:} At fulfillment of the conditions (\ref{dolgopoloveq:22}) the problem $C$ for
the equation (\ref{dolgopoloveq:mz1}) has unique solution represented by the formula (\ref{dolgopoloveq:21}).

\section{Introduction of the special representation\\ of the solution of the problem $\boldsymbol C$}

Let us introduce the special representation of the solution $W_p$ of the problem $C$ similar to what has been made by I.L.\,Karol~\cite{dolgopolovkarol} for the
Euler--Darboux equation on a plane.

For the purpose of convenience of our further reasoning, let us convert the formula of the solution
of the problem $C$ (\ref{dolgopoloveq:21}): rename the variables $x_0=x$, $y_0=y$,
$z_0=z$, $s=s_1$; let~us alternate the integration sequence in both integrals and do the replacement $s_1=t-s$, after that alternate the integrating sequence again. The solution of the problem~$C$ will come to the form
\begin{equation}
\begin{array}{c}
U(x,y,z)=\tau(x,y) -
\displaystyle\frac{1}{2}\displaystyle\frac{\Gamma(1-2p)}{\Gamma^2(1-p)}
\displaystyle\int\limits_{z}^{x-y}ds\displaystyle\int\limits_{y+s}^{x}\mu(t,t-s)\times\qquad\qquad\qquad\\
\qquad\qquad\times (x-t)^{-p}(t-y-s)^{-p} \,{
}_0F_1(1-p;-\lambda(x-t)(t-y-s))dt-\\
-\displaystyle\frac{\Gamma(2p)}{\Gamma^2(p)}
\displaystyle\int\limits_{z}^{x-y}(x-y-s)^{1-2p}ds\displaystyle\int\limits_{y+s}^{x}\nu(t,t-s)
(x-t)^{p-1}(t-y-s)^{p-1} \times\\
\qquad\qquad\qquad\qquad\qquad\qquad\qquad\qquad\times \,{
}_0F_1(p;-\lambda(x-t)(t-y-s))dt.
\end{array}
\label{dolgopoloveq:23}
\end{equation}

To simplify the formula (\ref{dolgopoloveq:23}), let us demand from the given function $\nu$ of integral~\cite{dolgopolovizv1} representation:
\begin{equation}
\nu(t,t-s) =
\displaystyle\int\limits_{0}^{t-s}T(s,\xi)(t-\xi-s)^{-2p} \,{
}_0F_1(1-p,\lambda(t-s-\xi)^2)d\xi.
\label{dolgopoloveq:24}
\end{equation}

{\bf Definition:} It is said that the solution
(\ref{dolgopoloveq:23}) of the problem $C$ for the equation
(\ref{dolgopoloveq:mz1})
has special representation
$W_p$ if
the function $\nu$
in it is defined by the equation
(\ref{dolgopoloveq:24}) where $T(x,y)$ is a new function,
$T(x,y)\in C^{(2)}(\overline D)$.

To find the form of the solution
$W_p$, substitute the function
(\ref{dolgopoloveq:24}) into the expression~(\ref{dolgopoloveq:23}) and
transform the third term. Let us alternate the integrating sequence in it over $t$ and~$\xi$
$$
\displaystyle\int\limits_{y+s}^xdt\displaystyle\int\limits_{0}^{t-s}d\xi
=
\displaystyle\int\limits_{0}^yd\xi\displaystyle\int\limits_{y+s}^{x}dt
+\displaystyle\int\limits_{y}^{x-s}d\xi\displaystyle\int\limits_{\xi+s}^{x}dt
$$
and consider the internal integrals
$$
\begin{array}{c}
J_1 =
\displaystyle\int\limits_{y+s}^x(x-t)^{p-1}(t-y-s)^{p-1}(t-\xi-s)^{-2p}\times\qquad\qquad\qquad\qquad\qquad\qquad\\
\qquad\qquad\qquad\qquad\times  \,{ }_0F_1(p,-\lambda(x-t)(t-y-s))
\,{ }_0F_1(1-p;\lambda(t-\xi-s)^2)dt.
\end{array}
$$

Let us represent the functions $\,{ }_0F_1(\alpha,z) =
\displaystyle\sum_{n=0}\displaystyle\frac{z^n}{(\alpha)_nn!}$
with the series and aply the rule on multiplication of the series, alternate the sequence
of summation and integration
$$
J_1 =
\displaystyle\sum_{k=0}^\infty\displaystyle\sum_{m=0}^k\displaystyle\frac{(-\lambda)^k(-1)^m}{m!(k-m)!(1-p)_m(p)_{k-m}}\displaystyle\int\limits_{y+s}^x(x-t)^{p-1+k-m}
\times\qquad\qquad\qquad\qquad
$$
$$
\qquad\qquad\qquad\qquad\qquad\times (t-s-y)^{p-1+k-m}(t-\xi-s)^{-2p+2m} dt.
$$
The integral will be named as $i$; using the transformation $t=x-(x-y-s)\mu$, express it via hyper-geometric Gauss function
$$
\begin{array}{c}
i=\displaystyle\frac{\Gamma^2(p)}{\Gamma(2p)}\displaystyle\frac{[(p)_{k-m}]^2}{(2p)_{2k-2m}}(x-y-s)^{2p-1+2k-2m}
\times\qquad\qquad\qquad\qquad\qquad\qquad\qquad\qquad\qquad\\
\qquad\qquad\times (x-\xi-s)^{-2p+2m}
F\left(2p-2m,p+k-m;2p+2k-2m;\displaystyle\frac{x-y-s}{x-\xi-s}\right).
\end{array}
$$

Name $\sigma=\displaystyle\frac{x-y-s}{x-\xi-s}$ and apply autotransform formula~\cite{dolgopolovBE73}
$$
F(\alpha,\beta,\gamma,\sigma)=(1-\sigma)^{\gamma-\alpha-\beta}F(\gamma-\alpha,\gamma-\beta;\gamma;\sigma).
$$
As a result, we have
$$
\begin{array}{c}
\!\!\!\!\!\!J_1 =
\displaystyle\frac{\Gamma^2(p)}{\Gamma(2p)}(x-y-s)^{2p-1}(x-\xi-s)^{-2p}\displaystyle\sum_{k=0}^\infty\displaystyle\frac{(-\lambda)^k}{(1-p)_kk!}\times
\qquad\qquad\qquad\qquad\qquad\qquad\\
\qquad\qquad\times (x-y-s)^{2k}(1-\sigma)^{k-p}
\displaystyle\sum_{m=0}^k\displaystyle\frac{(-1)^mk!(p)_{k-m}(1-p)_k}{m!(1-p)_m(k-m)!}\times\\
\null\hfill\times\displaystyle\frac{\sigma^{-2m}(1-\sigma)^m}{(2p)_{2k-2m}}F(2k;p+k-m;2p+2k-2m;\sigma).
\end{array}
$$

Using mathematical induction, we prove that the internal finite sum does not depend on $p$ and equals $\displaystyle\sum_{m=0}^k =
(-1)^k\sigma^{-2k}$, then
\begin{equation}
\begin{array}{c}
J_1=\displaystyle\frac{\Gamma^2(p)}{\Gamma(2p)}(x-y-s)^{2p-1}(x-\xi-s)^{-p}(y-\xi)^{-p}\times\qquad\qquad\qquad\\
\qquad\qquad\qquad\qquad\qquad\qquad\qquad\qquad\times  \,{
}_0F_1(1-p,\lambda(y-\xi)(x-\xi-s)).
\end{array}
\label{dolgopoloveq:25}
\end{equation}

Analogous reasoning makes it possible to obtain
\begin{equation}
\begin{array}{c}
J_2 =
\displaystyle\int\limits_{\xi+s}^x(x-t)^{p-1}(t-y-s)^{p-1}(t-\xi-s)^{-2p}
\,{ }_0F_1(p,-\lambda(x-t)(t-y-s))\times\\
\times \,{ }_0F_1(1-p,\lambda(t-\xi-s)^2)dt =
\displaystyle\frac{\Gamma(1-2p)\Gamma(p)}{\Gamma(1-p)}(x-s-y)^{2p-1}(\xi-s)^{-p}\times\qquad\\
\qquad\qquad\qquad\qquad\qquad\qquad\times  \,(x-\xi-s)^{-p}{
}_0F_1(1-p,-\lambda(x-\xi-s)(\xi-s)).
\end{array}
\label{dolgopoloveq:26}
\end{equation}

Recalling (\ref{dolgopoloveq:25}), (\ref{dolgopoloveq:26}) into the third
term of eq. (\ref{dolgopoloveq:23}) and substituting $t-s=\xi$ in the 2nd term containing
$\mu(t,t-s)$ follows
\begin{equation}
\begin{array}{l}
\!\!\!\!\!U(x,y,z)=\\
\!\!\!\!\!=\tau(x,y)-\!\!\displaystyle\int\limits_z^{x-y}\!\!\!ds\!\!\!
\displaystyle\int\limits_{0}^{y}\!\!\!T(s,\xi)(y-\xi)^{-p}(x-\xi-s)^{-p}\,{
}_0F_1(1-p,\lambda(y-\xi)(x-\xi-s))d\xi-\\
\null\hfill-\displaystyle\int\limits_z^{x-y}ds\displaystyle\int\limits_y^{x-s}N(s,\xi)(x-\xi-s)^{-p}(\xi-y)^{-p}\,{
}_0F_1(1-p,\lambda(y-\xi)(x-\xi-s))d\xi,
\end{array}
\label{dolgopoloveq:27}
\end{equation}
where
\begin{equation}
N(s,\xi)=\displaystyle\frac{1}{2\cos\pi p}
T(s,\xi)\, {\bf+}\, \displaystyle\frac{\Gamma(1-2p)}{2\Gamma^2(1-p)}\mu(\xi+s,\xi).
\label{dolgopoloveq:28}
\end{equation}

We need to note once again that the solution of the problem $C$ defined by the equation~(\ref{dolgopoloveq:21}) is obtained for the positive values of parameter $p$ while
for the negative values the sense is lost. The formula (\ref{dolgopoloveq:28}) is valid also
for $p=-q$. Really, with a direct check it is possible to show that the function
\begin{equation}
\begin{array}{c}
U(x,y,z) =
\tau_1(x,y)-\displaystyle\int\limits_z^{x-y}ds\displaystyle\int\limits_0^yT_1(s,\xi)(y-\xi)^q(x-\xi-s)^q\times
\qquad\qquad\qquad\qquad\\
\times  \,{ }_0F_1(1+q,\lambda(y-\xi)(x-\xi-s))d\xi-
\displaystyle\int\limits_z^{x-y}ds\displaystyle\int\limits_y^{x-s}N_1(s,\xi)(x-\xi-s)^q
\times\\
\quad\qquad\qquad\qquad\qquad\qquad\times (\xi-y)^q \,{
}_0F_1(1+q,\lambda(y-\xi)(x-\xi-s))d\xi,
\end{array}
\label{dolgopoloveq:29}
\end{equation}
where
\begin{equation} N_1(s,\xi)= \displaystyle\frac{1}{2\cos\pi q}
T_1(s,\xi)+\displaystyle\frac{\Gamma(1+2q)}{2\Gamma^2(1+q)}\mu_1(\xi+s,\xi)
\label{dolgopoloveq:30}
\end{equation}
is the solution of the equation (\ref{dolgopoloveq:mz1}) at $p=-q$
$\left(0<q<\displaystyle\frac{1}{2}\right)$, which satisfies
the conditions
\begin{equation}
U(x,y,x-y) =
\tau_1(x,y), \quad 0\leq y\leq x<+\infty;
\end{equation}
\begin{equation}
\lim_{z\to x-y-0}\displaystyle\frac{\partial U}{\partial z}=
\nu_1(x,y), \quad 0< y< x<+\infty;
\end{equation}
\begin{equation}
\lim_{z\to x-y-0} (x-y-z)^{-2q}(U_{xz}-U_{yz})=
\mu_1(x,y), \quad 0< y< x<+\infty,
\end{equation}
if the representation
\begin{equation}
\nu_1(t,t-s) =
\displaystyle\int\limits_0^{t-s}T_1(s,\xi)(t-\xi-s)^{2q} \,{
}_0F_1(1+q,\lambda(t-s-\xi)^2)d\xi,
\end{equation}
takes place where
$$
{\tau_1}_{xy}'' \in C(\overline D), \quad \nu_1, \mu_1, T_1 \in C(\overline D).
$$

As already stated, in the case of negative parameter Riemann method
proves to be inefficient, and even at $\lambda=0$, when there is
a possibility of obtaining the solution of Cauchy problem from the general solution
of the equation, it has complicated structure which in essence has lead to the necessity
of introduction of special
representations~\cite{dolgopolovkarol}.

Using the formulae (\ref{dolgopoloveq:27}), (\ref{dolgopoloveq:29}), it is possible to obtain
the solutions of new boundary value problems: Cauchy--Goursat, Darboux, problems with
integral conditions, and also problems
with pasting together
on
a singularity plane of the equation coefficients or on
characteristic plane.

\section{Solution of the mixed problem\\ in the special
representation~$\boldsymbol W_{\boldsymbol p}$}

As an example, let us consider the solution of one mixed problem for
the equation (\ref{dolgopoloveq:mz1}).

\textbf{Problem $C-G$:} (Cauchy--Goursat) in the region $\cal H$ to find
the solution of the equation~(\ref{dolgopoloveq:mz1}) continuous in $\overline{\cal H}$,
belonging to the class $W_p$, satisfying boundary conditions
(\ref{dolgopoloveq:2}), (\ref{dolgopoloveq:3}) and the condition
\begin{equation}
U(x,0,z)=\varphi(x,z),
\qquad 0\leq z\leq x<+\infty
\label{dolgopoloveq:35}
\end{equation}

Let us suppose that at $0<p<\displaystyle\frac{1}{2}$
$$
\begin{array}{c}\tau_{x,y} \in C(\overline D), \qquad \nu_{x,y} \in C^{(2)}(\overline D);\\ \nonumber
\tau(x,0)=0, \quad \varphi_{xz}'' \in C(\overline D), \quad \varphi(x,x)= \varphi_z'(z,z)=0. \end{array}
\eqno{(A)}
$$

Recall the equation (\ref{dolgopoloveq:27}) that is the solution
of the problem $C$
for eq. (\ref{dolgopoloveq:mz1}). This solution satisfies the conditions
(\ref{dolgopoloveq:2}), (\ref{dolgopoloveq:3}). We shall find the unknown function $N(s$,~$\xi)$
by assuming in the equation (\ref{dolgopoloveq:27}) $y=0$ and subordinating it
the condition (\ref{dolgopoloveq:35})
\begin{equation}
-\displaystyle\int\limits_z^xds\displaystyle\int\limits_0^{x-s}N(s,\xi)\xi^{-p}(x-\xi-s)^{-p}
\,{ }_0F_1(1-p,-\lambda(x-\xi-s)\xi)d\xi =\varphi(x,z).
\label{dolgopoloveq:36}
\end{equation}
The problem $C-G$ has reduced to integral equation in
$N(s,\xi)$. Supposing that the solution of the equation (\ref{dolgopoloveq:36})
exists, let us differentiate both sides of the identity (\ref{dolgopoloveq:36}) by~$z$, and then apply operator
$K_{p,x,z}[u]=\displaystyle\int\limits_z^{x+z}u(t,z)\,(x+z-t)^{p-1}dt$ to both sides of the obtained equation preliminarily
replacing $x$ with $t$ in eq. (\ref{dolgopoloveq:36}).

As a result of calculations,obtain
\begin{equation}
\displaystyle\int\limits_0^xN_1(z,\xi) \,{
}_0F_1(1,-\lambda\xi(x-\xi))d\xi = f(x,z),
\label{dolgopoloveq:37}
\end{equation}
where
\begin{equation}
N_1 = \xi^{-p}N, \qquad f(x,z) =
\displaystyle\frac{1}{\Gamma(p)\Gamma(1-p)}\displaystyle\int\limits_z^{x+z}\varphi_z'(t,z)(x+z-t)^{p-1}dt.
\label{dolgopoloveq:38}
\end{equation}
The unique solution of the equation (\ref{dolgopoloveq:37})
has been obtained in the paper~\cite{dolgopolovBR07}
\begin{equation}
N_1(z,x) = f_x'(x,z) -
x^{-1}\displaystyle\int\limits_0^xf_t'(t,z)t\displaystyle\frac{\partial}{\partial
t} \,{ }_0F_1(1,\lambda x(x-t))dt
\end{equation}
on condition $f(0,z)=0$.

Realling the expression of the function $f$ from the formula (\ref{dolgopoloveq:38}) into the formula (\ref{dolgopoloveq:37}) and performing the calculations, we have
\begin{equation}
\begin{array}{c}
N(z,x)=\displaystyle\frac{x^p}{\Gamma(p)\Gamma(1-p)}\left[\displaystyle\int\limits_z^{x+z}\varphi_{zt}''(t,z)(x+z-t)^{p-1}
\,{ }_0F_1(p,\lambda x(x-t+z))dt-\right.\\
\left. \qquad\qquad\qquad\qquad\qquad-
\displaystyle\frac{\lambda}{p}\displaystyle\int\limits_z^{x+z}\varphi_{z}'(t,z)(x+z-t)^{p}
\,{ }_0F_1(1+p,\lambda x(x-t+z))dt\right].
\end{array}
\label{dolgopoloveq:40}
\end{equation}

%%%%%%%%%%%%%%%%%%%%%%%%%%%%%%%%%%%%%%%%%%%%%%%%%%%%%%%%%%%%%%%%%%%%%%%%%%%%

Direct checking shows that at satisfying the conditions $(A)$, the function (\ref{dolgopoloveq:40}) is the unique solution of eq.~(\ref{dolgopoloveq:37}).

In case of negative parameter of the equation (\ref{dolgopoloveq:mz1}) $p=-q$, the solution of the problem $C-\cal G$ is reduced to integral equation
in $N_1$:
\begin{equation}
-\displaystyle\int\limits_z^{x-y}ds\displaystyle\int\limits_y^{x-s}N_1(s,\xi)(x-\xi-s)^q\xi^q
\,{ }_0F_1(1+q,-\lambda (\xi-s)(x-\xi-s))d\xi = \varphi_1(x,z).
\label{dolgopoloveq:43}
\end{equation}

Supposing that a solution of the equation (\ref{dolgopoloveq:43}) exists,
let us differentiate the identity~(\ref{dolgopoloveq:43}) first with respect to $z$, then to
$x$ and apply operator $K_{q,x,z}[u]$ to both parts of the obtained equation. Then,
reasoning in a similar way as in the case of eq. (\ref{dolgopoloveq:36}), it is possible to prove that
the function
\begin{equation}
\begin{array}{c}
N_1(z,x)=\displaystyle\frac{x^{-q}}{\Gamma(1+q)\Gamma(1-q)}\left[\displaystyle\int\limits_z^{x+z}{\varphi_1}_{ztt}'''(t,z)(x+z-t)^{-q}
\,{ }_0F_1(1-q,\lambda x(x-t+z))dt-\right.\\
\left.
\qquad\qquad\qquad-\displaystyle\frac{\lambda}{1-q}\displaystyle\int\limits_z^{x+z}{\varphi_1}_{zt}''(t,z)(x+z-t)^{1-q}
\,{ }_0F_1(2-q,\lambda x(x-t+z))dt\right]
\end{array}
\label{dolgopoloveq:44}
\end{equation}
is the unique solution of the eq. (\ref{dolgopoloveq:43}) at complimentary
to $(A)$ conditions superimposed on the given function $\varphi_1$ which have the form:
\begin{equation}
{\varphi_1}_{xzx}'' \in C(\overline D), \qquad
{\varphi_1}_{zx}''(z,z) = 0.
\label{dolgopoloveq:45}
\end{equation}

{\bf Theorem 2:} The problem $C-\cal G$ for eq. (\ref{dolgopoloveq:mz1})
has a unique solution defined by the equations (\ref{dolgopoloveq:27}),
(\ref{dolgopoloveq:40}) at fulfilment of conditions $(A)$ in case
$0<p<\displaystyle\frac{1}{2}$ and defined by the formulae
(\ref{dolgopoloveq:29}), (\ref{dolgopoloveq:44}) at fulfilment of conditions $(A)$ and
(\ref{dolgopoloveq:45}) for $p=-q$
$\left(0<q<\displaystyle\frac{1}{2}\right)$.

\def\refname{
%%\Large
References}

\end{document}